\documentclass{ifacconf}

\usepackage{natbib}            
\usepackage{graphicx}          
\usepackage{psfrag}
\usepackage{amsmath}
\usepackage{multirow}
\usepackage[latin1]{inputenc}
\usepackage[T1]{fontenc}
\usepackage{amssymb,amsfonts,textcomp}
\usepackage{url}

\usepackage{pstricks}
\usepackage{psfrag}                               
\usepackage[dvips]{epsfig}    

\newcounter{definition}

\newcounter{remark}

\newcounter{example}

\newenvironment{ex*}[1]{\addtocounter{example}{1}\smallbreak\noindent
  {\bf Example \theexample{} --- {\bf #1}}}{\hfill$\Box$\smallbreak}
\newcommand{\realR}{\mathbb{R}}

\newcommand{\RL}{{\bf RL}_{\infty}}

\newcommand{\RHii}[2]{{\bf RH}_{\infty}^{{#1}\times{#2}}}
\newcommand{\RLii}[2]{{\bf RL}_{\infty}^{{#1}\times{#2}}}

\DeclareMathOperator{\imag}{Im}

\DeclareMathOperator{\RE}{Re}
\begin{document}

\begin{frontmatter}
\title{\text{\!\!\!\!\!\!\!\!\!\!\!\!\!\! Structure Preserving H-infinity Optimal PI Control}} 


\author[First]{Anders Rantzer} 
\author[Second]{Carolina Lidström} 
\author[Third]{Richard Pates}
\address[First]{Lund University, Sweden (e-mail: rantzer@control.lth.se).}                                              \address[Second]{Lund University, Sweden (e-mail: carolina.lidstrom@control.lth.se).}
\address[Third]{Lund University, Sweden (e-mail: richard.pates@control.lth.se).}

\begin{abstract}                          
A multi-variable PI (proportional integrating) controller is proved to
be optimal for an important class of control problems where
performance is specified in terms of frequency weighted H-infinity norms. 
The problem class includes networked systems with a subsystem in each
node and control action along each edge. For such systems, the optimal
PI controller is decentralized in the sense that control
action along a given network edge is entirely determined by states at nodes
connected by that edge. 
\end{abstract}
\end{frontmatter}
\section{Introduction}
Classical theory for multi-variable control synthesis suffers from a
severe lack of scalability. Not only does the computational cost for
Riccati equations and LMIs grow rapidly with the state dimension, but
also implementation of the resulting controllers becomes unmanageable
in large networks due to requirements for communication, computation
an memory. Because of this situation, considerable research efforts
have recently been devoted to development of scalable and structure preserving
methods for design and implementation of networked control
systems building on early contributions by \cite{BamiehPaganini+02}, \cite{d2003distributed} and \cite{Rotkowitz+02}.

In practice, most scalable control architectures are built on layering.
For example,
control systems in the process industry are often organized in 
a hierarchical manner, where scalar PID controllers are used at the lowest
level and the reference values of these controllers are computed on a
slower time scale by centralized optimization algorithms. This approach
often works well provided that the scalar loops are reasonably
decoupled and that the coordination dynamics are comparatively
slow. Other important applications where decentralized and layered control
architectures have been successfully applied are power systems and
Internet traffic control.

A scalable approach to control synthesis has recently been developed
based on the notion of 
positive systems (\cite{rantzer2015ejc}) and the nonlinear
counterpart monotone systems (\cite{7753319}). Such systems are characterized by
existence of Lyapunov functions and other performance certificates
whose complexity grows only linearly with the systems size.
Interestingly, the maximal gain is always attained at zero
frequency. Restricting attention to closed loop positive systems, 
distributed static controllers were optimized subject to
$H_\infty$ performance  in \cite{Tanaka+11}, while $\ell_1$
performance was considered in \cite{Briat13}.

The powerful synthesis methods for positive systems have raised an important
question: How restrictive is a demand for closed loop positivity in
$H_\infty$ optimal control?
It was therefore a remarkable step forward when \cite{lidstrom+16acc}
showed that for a large class of networked control systems with
``diffusive'' dynamics (i.e. symmetric state matrix), it is not
restrictive at all. Instead, controllers defined by a simple
closed form expression achieve
the same level of performance as centralized controllers derived using
Riccati equations or LMIs. In particular, the following problem was
considered:

\emph{Given a graph $(\mathcal{V},\mathcal{E})$ and the system 
\begin{align}
  \notag\\[-4mm]\dot{x}_i&=a_ix_i+\sum_{(i,j)\in\mathcal E}(u_{ij}-u_{ji})+w_i&i&\in\mathcal{V},
\label{eqn:netsys}
\end{align}
find a control law of the form $u=Kx$ that minimizes the
$H_\infty$ norm of the transfer function from the disturbance $w$ to
the controlled output $(x,u)$.} 

\cite{lidstrom+16acc} proved that when $a_i<0$ an optimal
control law is given by
\begin{align}
  \notag\\[-4mm] u_{ij}&=x_i/a_i-x_j/a_j \label{eqn:Pcontrol} \\[-4mm]\notag
\end{align}
and the closed loop from $w$ to $x$ is a positive
system. This control law is trivial to 
compute and decentralized in the sense that control
action on the edge $(i,j)$ is entirely determined by the states in
node $i$ and node $j$. 

The expressions above define a rare but important class of systems where decentralized controllers
are known to achieve the same $H_\infty$ performance as the best
centralized ones. Still, the setting is insufficient for many
practical applications. In particular, as all proportional
controllers, the control law (\ref{eqn:Pcontrol}) is unable to remove
static errors in presence of constant disturbances. The purpose of
this paper is to get rid of this deficiency, by modifying the performance
criterion to optimize dynamic controllers with integral action.

The structure of the paper is as follows: After introduction of some
basic notation in section 2, we prove the main result in
section 3,  consider some basic applications in section 4 and 5, before
making conclusions in section~6. A well known
matrix optimization result is included for completeness as a short appendix.

\section{Notation}
Let $\RL$ be the set of proper (bounded at
infinity) rational functions with real coefficients. 
The set of $m\times n$ matrices with elements in
$\RL$ is denoted $\RLii{m}{n}$. Given ${P}\in\RLii kl$, we say
that ${K}\in\RLii lk$ is \emph{stabilizing} ${P}$ provided that
\begin{align*}
  \begin{bmatrix}I\\{K}\end{bmatrix}
  \left[I+{P}{K} \right]^{-1}
  \begin{bmatrix}I&{P}\end{bmatrix}\;\;\in \;\; \RHii{(k+l)}{(k+l)}
\end{align*}
has no poles in the closed right half plane. Furthermore
$\|P\|_\infty=\sup_{\RE s\ge0}\|P(s)\|$. For a matrix
$M\in\realR^{n\times m}$, we denote the pseudo-inverse by $M^\dagger$
and the spectral norm by $\|M\|$. For a square symmetric matrix $M$
the notation $M\succ0$ means that $M$ is positive definite, while
$M\prec0$ means that $M$ is negative definite.

\section{Main result}

\begin{thm}
  Let ${P}(s)=(sI-A)^{-1}B$ with $A$ 
  symmetric negative definite. Assume that $\tau\ge\sqrt{\|B^TA^{-4}B\|}$.  
Then the problem
\begin{center}
  \begin{tabular}{ll}
    \\[-3mm]
    Minimize &$\|\left(I+KP\right)^{-1}{K}\|_\infty$\\[2mm]
    subject to & $\|\frac{1}{s}P(I+{K}{P})^{-1}\|_\infty\le\tau$\\[2mm]
  \end{tabular}
\end{center}
over stabilizing $K$, is solved by
\begin{align*}
  \widehat{{K}}(s)&=k\left(B^TA^{-2}-\frac{1}{s}B^TA^{-1}\right)
\end{align*}
where $k=\|(A^{-1}B)^\dagger\|/\tau$.
\label{thm:main}
\end{thm}

\begin{pf}
Define ${F}_{{K}}=\left(I+{K}{P}\right)^{-1}{K}$ and factorize $B$ as $B=GH^T$, where $G$ and $H$ have full column rank. Then
\begin{align*}
  {F}_{\widehat{{K}}}(s)
  &=k\left(s+kB^TA^{-2}B\right)^{-1}B^TA^{-2}(sI-A)\\
  &=kH\left(s+kG^TA^{-2}GH^TH\right)^{-1}G^TA^{-2}(sI-A)
\end{align*}
so the poles of $F_K$ are the eigenvalues of
$-kG^TA^{-2}GH^TH$, which are equal to the non-zero eigenvalues of $kB^TA^{-2}B$.
In particular, $F_{\widehat{{K}}}$ and $P$ are stable, so ${\widehat{{K}}}$ is stabilizing.

We also have
\begin{align*}
  I&\preceq \tau^2k^2B^TA^{-2}B\\
  B^T(\omega^2I+A^2)^{-1}B&\preceq\tau^2\left[\omega^2I+k^2(B^TA^{-2}B)^2\right]
\end{align*}
so
\begin{align*}
  \tau&\ge\left\|(i\omega I-A)^{-1}B(i\omega I+kB^TA^{-2}B)^{-1}\right\|\\
  &=\left\|\frac{1}{i\omega}P(i\omega)(I+\widehat{K}{(i\omega)}{P(i\omega)})^{-1}\right\|.
\end{align*}
In general
\begin{align*}
  PF_KP=P-P(I+KP)^{-1}
\end{align*} 
so the constraint $\|\frac{1}{s}P(I+{K}{P})^{-1}\|_\infty\le\tau$ gives
\begin{align*}
  P(0)F_K(0)P(0)=P(0).
\end{align*}
Hence consider the minimization problem at $s=0$, i.e. to minimize
  $\|F\|$ subject to $P(0)FP(0)=P(0)$. Standard
  calculations (Lemma~\ref{lem:app} in the appendix) gives that the minimal value 
  $\gamma={\|P(0)^\dagger\|}$ is attained by
  $F=P(0)^\dagger$. In particular, since
  ${F}_{\widehat{{K}}}(0)={P}(0)^\dagger$, it follows that $\widehat{K}$ is
  optimal at $\omega=0$.
 
The inequality $\|{F}_{\widehat{{K}}}(i\omega)\|\le
\|{F}_{\widehat{{K}}}(0)\|$ can be rewritten as
\begin{eqnarray*}
&{F}_{\widehat{{K}}}(i\omega)^*{F}_{\widehat{{K}}}(i\omega)
  \preceq \gamma^2 I\\
&A^{-2}B\left[\omega^2k^{-2}I+(B^TA^{-2}B)^2\right]^{-1}B^TA^{-2}
  \preceq \gamma^2(\omega^2I+A^2)^{-1}\\
&k^2B^TA^{-2}(\omega^2I+A^2)A^{-2}B\preceq \gamma^2\left[\omega^2I+k^2(B^TA^{-2}B)^2\right]
\end{eqnarray*}
The last inequality holds trivially for $\omega=0$ and it holds
for all other $\omega$ provided that $k^2B^TA^{-4}B \preceq \gamma^2I$,
which is equivalent to the assumption on $\tau$. Thus
$\|{F}_{\widehat{{K}}}\|_\infty$ takes the minimal value
$\gamma$ and the proof is complete.
\end{pf}

\section{Control on Networks}

Theorem~\ref{thm:main} can be applied to the following problem:

\emph{Given a graph $(\mathcal{V},\mathcal{E})$, suppose that
\begin{align}
  \left\{
  \begin{array}{rl}
    \dot {x}_i&=\;a_ix_i+b_iu_i+\hbox{$\sum_{(i,j)\in\mathcal E}$}(u_{ij}+v_{ji})\\
    e_i&=\;r_i-x_i
  \end{array}\right.
\label{eqn:netsys2}
\end{align}
$i\in\mathcal{V}$, $a_i<0$, $x(0)=0$ and $u_{ij}=-u_{ji}$,
$v_{ij}=-v_{ji}$. The signals $u$, $v$, $r$ and $e$ can be viewed as input-,
\hbox{disturbance-}, reference- and error-signals respectively. 
The problem is to find a control law $u=k*e$ that minimizes the
$L_2$-gain from $r$ to
$u$ while keeping the $L_2$-gain from $v$ to $z$
bounded by $\tau$ when $\dot{z}=x$, $z(0)=0$.} 

The optimal controller $\widehat{K}$ has a distributed realization with
one integrator at every node:
\begin{align}
  \left\{
  \begin{array}{rl}
    \dot{z}_i &=\;e_i\\
    u_{ij}&=\;z_i/a_i-e_i/a_i^2-z_j/a_j+e_j/a_j^2\\
    u_i&=\;b_i(z_i/a_i-e_i/a_i^2) 
  \end{array}\right.
  \label{eqn:decentPI}
\end{align}
Just like (\ref{eqn:Pcontrol}), this control law is decentralized in the sense that control
action on the edge $(i,j)$ is entirely determined by the errors at
nodes $i$ and $j$. The closed loop map from $r$ to $z$
has transfer matrix
$$(sI-A)^{-1}(sI+kBB^TA^{-2})^{-1},$$ and non-negative impulse
response. It should be noted that unless the controller realization is
minimal, the controller will have integrators that are not stabilized
in closed loop. For this reason, it necessary that $b_i\ne0$ for at least
one node in every connected component of the graph.

\bigbreak

\section{Example}
\begin{figure} 
\centering
\includegraphics[width=.7\hsize]{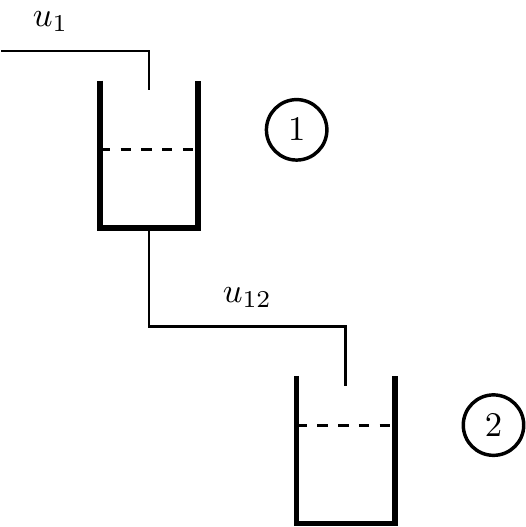}
\caption{System with two buffers 1 and 2. The interconnection given
  by $B$ in \eqref{matrices} is depicted by the line drawn between the
  buffers as well as between buffer 1 and the
  exterior. Degradation due to internal dynamics is not
 shown. \vskip 5mm }\label{fig:bufferSys}
\end{figure}


Consider the system depicted in Figure~\ref{fig:bufferSys}. The dynamics of the levels in the two buffers marked 1 and 2, around some steady state, is given by 
\begin{equation}
\begin{aligned}
\dot{x}_1 &= -x_1 +u_1-u_{12}&e_1&=r_1-x_1\\ 
\dot{x}_2 &= -2x_2 +u_{12} &e_2&=r_2-x_2,
\end{aligned} \label{systemBuff}
\end{equation}
where $x_1$ is the level in buffer 1 and $x_2$ is the level in buffer 2. The transfer function of \eqref{systemBuff} is given by ${P(s) = (sI-A)^{-1}B}$ with 
\begin{align}
A &= \begin{pmatrix}-1 & 0 \\ 0 & -2 \end{pmatrix}& B &= \begin{pmatrix} 1 & -1 \\ 0 & 1 \end{pmatrix}. \label{matrices}
\end{align}
The criterion on $\tau$ is given by $\tau \geq 1.43$ and $\widehat{K}$ is given by
\begin{align*}
u_1(t) &= -ke_1  -k\int_0^t e_1(\sigma) \, d \sigma \\ 
u_{12}(t) &= ke_1-\frac{k}{2}e_2 + k\int_0^t\left(e_1(\sigma)-\frac{1}{4}e_2 (\sigma)\right)\, d\sigma. \\  
\end{align*}
for $k = 8.5/\tau$. Notice that each control input is only using the
state(s) it affects through the matrix $B$, i.e., the controller has
the same zero-block structure as $B^T$.
Thus, the controller only considers local information, where local is subject to the interconnection specified by $B$. 

\begin{figure} 
\begin{center}
\includegraphics[width=.8\hsize]{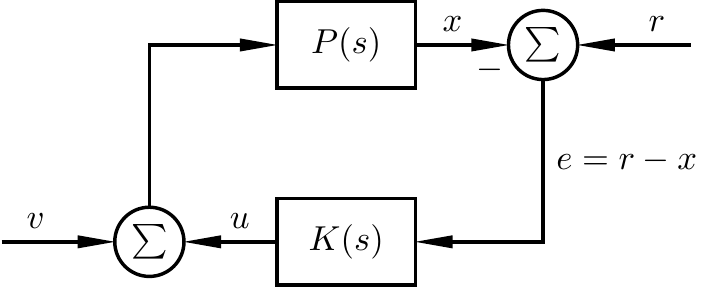}
\end{center}
\caption{Block diagram of closed-loop system with
  measurement $y$, control input $u$, process disturbance $v$,
  reference  value (or measurement error) $r$, controller~$K$ and process $P$.}\label{blockdiagram}
\end{figure}

The optimization objective stated in Theorem~1 concerns the transfer
function $$(I+KP)^{-1}K,$$ which maps the
reference value $w$ to the control input $u$, as depicted in
Figure~\ref{blockdiagram}. Thus, the first objective is to minimize the
control effort needed to follow the reference value $w$. The second performance criterion 
\begin{align}
\left \| \frac{1}{s}P(I+KP)^{-1}\right \|_{\infty} \leq \tau \label{critKP}
\end{align}
specifies the control quality in terms of disturbance rejection. The impact from a low frequency process disturbance $v$ should
be attenuated by the feedback loop. The parameter $\tau$ is a time
constant that determines the bandwidth of the control loop. The impact
of $\tau$ will be illustrated below. 

\begin{figure} 
\centering
\includegraphics[width=.9\hsize]{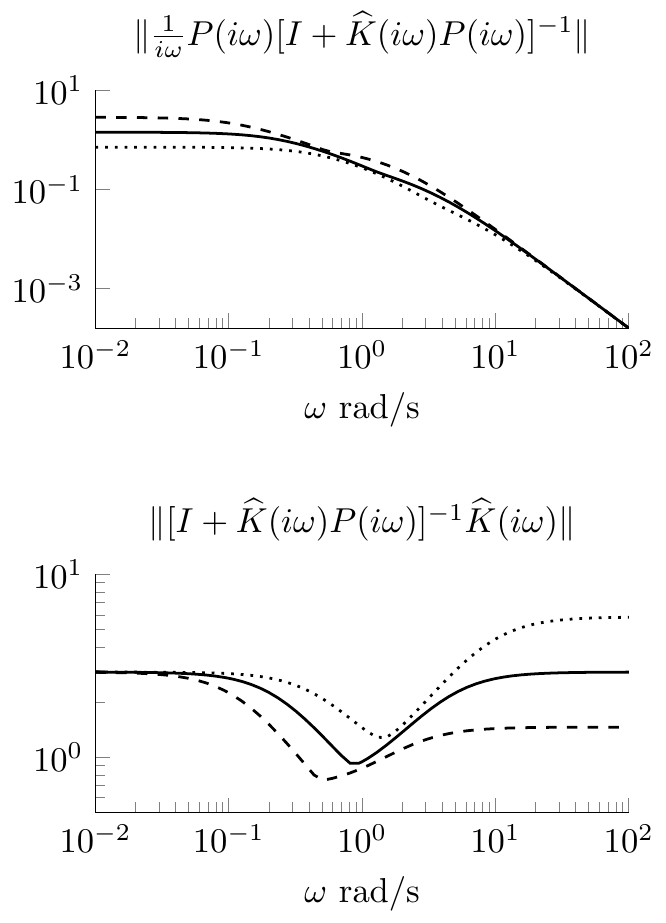} 
\bigskip
\centering
\includegraphics[width=\hsize]{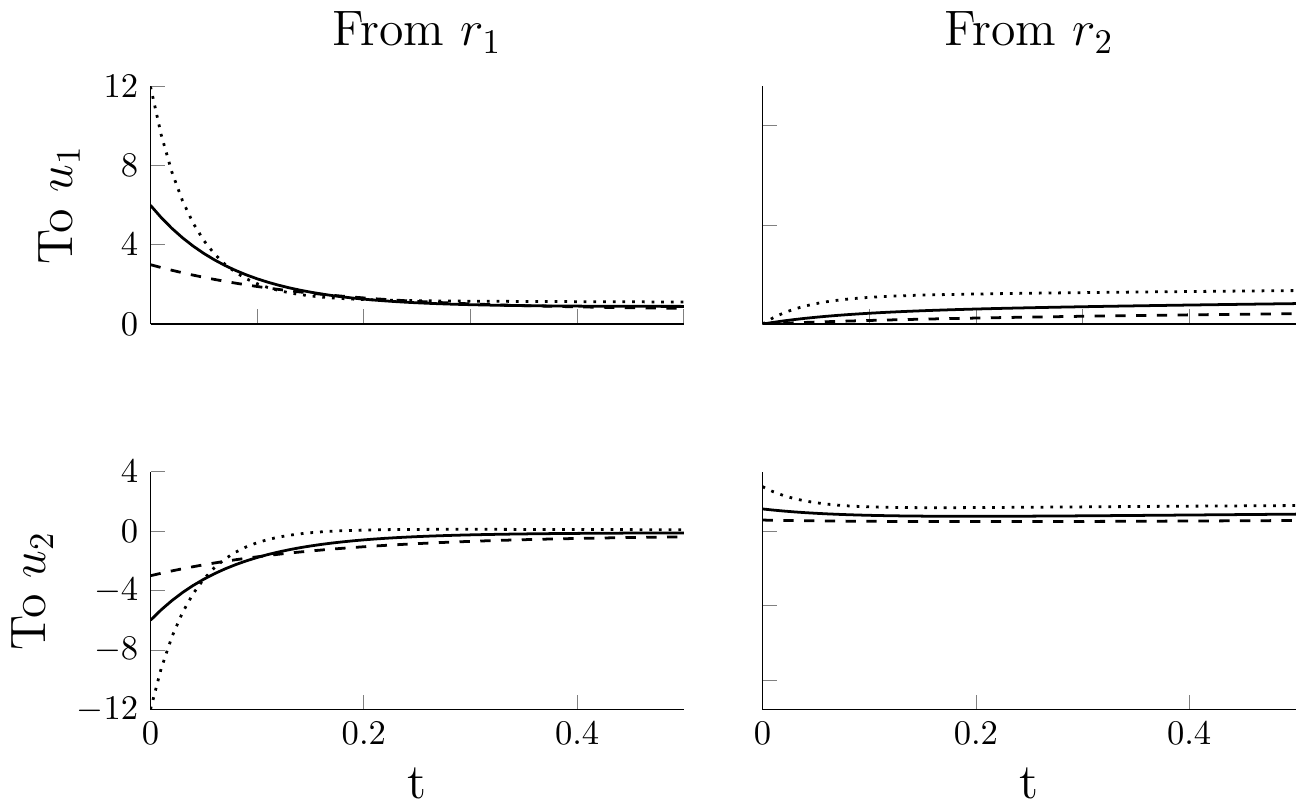} 
\caption{Spectral norms and step responses for different $\tau$; define $\tau_* =
  \sqrt{\|B^TA^{-4}B\|}$, then $\tau = \tau_*$ gives the solid line,
  $\tau = 2\tau_*$ the dashed line and $\tau = \tau_*/2$ the dotted
  line. The last case violates the optimality conditions of
  Theorem~\ref{thm:main} and the norm is not maximal at $\omega=0$.}
 \label{fig:differentTau}
\end{figure}

Given 
\begin{align*}
  P(s) &= (sI-A)^{-1}B\\
  \widehat{K}(s) &= \frac{\|(A^{-1}B)^\dagger\|}{\tau}\left(B^TA^{-2}-\frac{1}{s}B^TA^{-1} \right),
\end{align*}
Figure~\ref{fig:differentTau} plots the norms of
$\|\frac{1}{i\omega}P(i\omega)[I+\widehat{K}(i\omega)
P(i\omega)]^{-1}\|$ and $\|[I+\widehat{K}(i\omega)
P(i\omega)]^{-1}\widehat{K}(i\omega)\|$ against the frequency $\omega$
for three different values of
$\tau$. 

The first diagram clearly shows that the disturbance rejection
is increasingly effective as the time constant $\tau$ is
reduced. However, as shown in the second diagram, this comes at the
price of larger control signals at higher frequencies, while the gain at
$\omega=0$ remains unchanged. In the step response diagrams, a smaller
value of $\tau$ results in larger control values for small $t$, while
the steady state values remain unchanged.

\section{Conclusions}
This paper has formulated a class of dynamic state feedback control
problems for which a structure preserving PI controller is $H_\infty$ optimal.
An explicit expression
for the optimal gain has been given, which clarifies the
relationship between plant structure and achievable preformance.

\begin{ack}                               
This work was supported by the Swedish Research Council through the LCCC Linnaeus Center. The authors are members of the LCCC Linnaeus Center and the eLLIIT Excellence Center at Lund University. 
\end{ack}


\section{Appendix}

\label{App1}   
\begin{lem}
Let $A\in\mathbb{C}^{n\times{}m}$. Then
\begin{equation}\label{eq:matopt}
\begin{aligned}
\min_{X\in\mathbb{C}^{m\times{}n}}&\|X\|
\;\;\; \mathrm{s.t.}\;\;\; AXA=A
\end{aligned}
\end{equation}
has the minimal value $\|A^\dagger\|$, attained by $\widehat{X}=A^\dagger$.
\label{lem:app}
\end{lem}
\begin{pf}
Let $y\in\imag(A)$ be a unit vector in $\mathbb{C}^m$ and let $X$ be any feasible point for \eqref{eq:matopt}. Consider:
\begin{equation}\label{eq:matopt1}
\begin{aligned}
\min_{x\in\mathbb{C}^{n}}&|x|
\;\;\; \hbox{s.t.}\;\;\; y=Ax \hbox{ and }x=Xy.
\end{aligned}
\end{equation}
Observe that because $X$ is feasible , the optimal solution always
exists. The value equals $|Xy|$, giving a lower bound for $\|X\|$. Relaxing
(\ref{eq:matopt1}) by removing the second constraint gives a
least squares problem with solution $x=A^\dagger y$ and value
$|A^\dagger y|$. Maximizing over $y$ gives $\|A^\dagger\|$. The result
follows because $X=A^\dagger$ achieves this lower bound.
\end{pf}

\end{document}